\documentclass[10pt,a4paper]{article}
\usepackage{a4,amsfonts,amssymb,amscd,amsmath,bbm,amsthm,cite}
\usepackage{xypic,nicefrac,epsfig}

\theoremstyle{definition}

\newtheorem{Thm}{Theorem}[section]

\newtheorem{Lemma}[Thm]{Lemma}
\newtheorem{Cor}[Thm]{Corollary}

\newtheorem{Assumption}[Thm]{Assumption}
\newtheorem{Example}[Thm]{Example}

\newcommand{\bN}{\ensuremath{\mathbbm{N}}}
\newcommand{\bNp}{\ensuremath{\mathbbm{N}^+}}

\newcommand{\bR}{\ensuremath{\mathbbm{R}}}
\newcommand{\bRp}{\ensuremath{\mathbbm{R}^+}}
\newcommand{\bRpO}{\ensuremath{\mathbbm{R}^+_0}}

\newcommand{\cU}{\ensuremath{{\mathcal U}}}

\newcommand{\Lin}{\ensuremath{{\mathcal L}}}

\newcommand{\myfrac}[2]{\ensuremath{\nicefrac{#1}{#2}}}
\newcommand{\norm}[1]{{\ensuremath{\lvert#1\rvert}} }
\newcommand{\opnorm}[1]{{\ensuremath{\lVert#1\rVert}} }
\newcommand{\supnorm}[1]{{\ensuremath{{\lVert#1\rVert}_{\infty}}} }
\newcommand{\lunendl}{\ensuremath{{l^{\infty}}}}
\newcommand{\fnto}{\rightarrow}
\newcommand{\limto}{\rightarrow}
\newcommand{\la}{\langle}
\newcommand{\ra}{\rangle}
\newcommand{\fI}{\ensuremath{\langle f_i\rangle}}
\newcommand{\hI}{\ensuremath{\langle h_i\rangle}}

\newcommand{\xN}{\ensuremath{{\mathbbm{^{\ast}\bN}}}}
\newcommand{\xR}{\ensuremath{{\mathbbm{^{\ast}\bR}}}}
\newcommand{\xRp}{\ensuremath{{\x\bR^+}}}

\newcommand{\x}{\ensuremath{^{\ast}}}
\newcommand{\xT}{\ensuremath{{{^{\ast}T}}}}
\newcommand{\xD}{\ensuremath{{{^{\ast}D}}}}
\newcommand{\xA}{\ensuremath{{^{\ast}A}}}
\newcommand{\xDA}{\ensuremath{{^{\ast}D(A)}}}

\newcommand{\dach}[1]{\ensuremath{\widehat{#1}}}
\newcommand{\fD}{\ensuremath{\dach{f}}}
\newcommand{\gD}{\ensuremath{\dach{g}}}
\newcommand{\TD}{\ensuremath{\dach{T}}}
\newcommand{\TT}{\ensuremath{\tilde{T}}}
\newcommand{\AD}{\ensuremath{\dach{A}}}
\DeclareMathOperator{\fin}{fin}
\DeclareMathOperator{\std}{std}
\newcommand{\near}{\ensuremath{\approx}}

\newcommand{\MAX}[2][T]{\ensuremath{{\left(#2\right)}^{#1\text{--max}}}}

\newcommand{\xE}{\ensuremath{{^{\ast}E}}}
\newcommand{\finxE}{\ensuremath{{\fin{^{\ast}E}}}}
\newcommand{\EO}{\ensuremath{E_0}}
\newcommand{\ED}{\ensuremath{\widehat{E}}}
\newcommand{\ET}{\ensuremath{E_T}}
\newcommand{\ETD}{\ensuremath{\widehat{\ET}}}
\newcommand{\Emax}{\ensuremath{E_{\text{max}}}}
\newcommand{\EmaxD}{\ensuremath{\widehat{E}_\text{max}}}

\newcommand{\luE}{\ensuremath{\lunendl(E)}}
\newcommand{\mT}{\ensuremath{m^T}}
\newcommand{\CU}{\ensuremath{c_{\mathcal U}}}
\newcommand{\luED}{\ensuremath{\myfrac{\luE}{\CU}}}
\newcommand{\mTD}{\ensuremath{\myfrac{\mT}{\CU\cap\mT}}}

\newcommand{\myisomorphism}{\ensuremath{\iota}}  
\newcommand{\myisometry}{\ensuremath{\phi}}  

\newcommand{\mycc}{C^\text{b}} 

\begin{document}

\begin{center}
{\Large F--Products and Nonstandard Hulls for Semigroups}

{\large Jakob Kellner, Hans Ploss}
\end{center}

\section{Abstract}

Derndinger \cite{derndinger} and Krupa \cite{krupa} defined the
F--product of a (strongly continuous one--parameter) semigroup 
(of linear operators) and presented some applications (e.g. to spectral
theory of positive operators, cf. \cite{engelnagel}).  Wolff (in
\cite{wolff_spectral} and \cite{wolff_nsfa}) investigated some kind of
nonstandard analogon and applied it to spectral theory of group
representations. The question arises in which way these constructions are
related.

In this paper 
we show that the classical and the
nonstandard F--product are isomorphic 
(theorem \ref{main}).
We also prove a little ``classical''
corollary (\ref{little}).

\section{Basic Notation}\label{basicnot}

\subsection{Semigroups:}

Let $E$ be a (real or complex) Banachspace.
We denote the norm of an element $f$ of $E$ by 
$\norm{f}$. $\Lin(E)$ is the set of 
bounded linear functions from $E$ to $E$.
The elements of $\Lin(E)$ are called bounded linear 
{\em operators},
and the norm of an operator $A$ is 
is denoted by $\opnorm{A}$.
A one--parameter semigroup of bounded linear operators
(or {\em semigroup}, for short) is 
a function $T: \bRpO\fnto \Lin(E)$ such that 
$T(0)=\text{Id}_E$ and $T(t_1+t_2)=T(t_1)\circ T(t_2)$.
A semigroup $T$ is called strongly continuous
(or {\em continuous}, for short), if for all $f\in E$,
$\lim_{t\limto 0}\norm{T(t)f-f}=0$
(i.e. $\lim_{t\limto 0}T(t)=\text{Id}_E$ in the strong
topology).
$T$ is called {\em uniformly
continuous},
if $\lim_{t\limto 0}\opnorm{T(t)-\text{Id}_E}=0$
(i.e. $\lim_{t\limto 0}T(t)=\text{Id}_E$ in the topology 
induced by the operator--norm).

For a continuous semigroup $T$ and $f\in E$ we define 
$A(f):=\lim_{h\limto 0}\frac{1}{h}(T(h)f-f)$, if this limit exists
(in $E$). $D(A)$ is the set of $f\in E$ such that $A(f)$ exists.
$A$ is called the {\em Generator} of $T$.

%
More about the (classical) theory of continuous semigroups can be 
found in the textbook \cite{engelnagel}.
We will only need the following results:

\begin{Lemma}\label{lem:basicsemigroup}
\begin{enumerate}
\item\label{itembound} 
For every continuous semigroup $T$ there exist
constants $M\geq 1$, $\omega\in \bR$ such that for all $t\in\bRpO$:
$\opnorm{T(t)}\leq M e^{\omega t}$.
\item
$D(A)$ is a dense linear subspace of $E$, and $A$
is a closed linear operator (i.e. the graph of $A$ is
a closed subset of $E\times E$).
\item
$D(A)=E$ if and only if  $T(t)$ is uniformly continuous.
\item\label{itemXY1}
If $f\in D(A)$ and $h>0$, then
$\norm{T(h)f-f}\leq 
h \cdot \norm{A(f)}\cdot \sup_{s\leq h}\opnorm{T(s)}$
\item\label{itemXY2}
If $f\in D(A)$ and $h>0$, then
$\norm{\frac{T(h)f-f}{h}-A(f)}\leq %
\sup_{s\leq h}(\norm{T(s)A(f)-A(f)})$
\item\label{itemuoehgqtetqw}
For any $\bar{f}$ and $h>0$ there is a $f\in D(A)$ such that 
$A(f)=\frac{T(h)\bar{f}-\bar{f}}{h}$ and $\norm{\bar{f}- f}\leq 
 \sup_{s\leq h}(\norm{T(s)\bar{f}-\bar{f}})$

\end{enumerate}
\end{Lemma}

Remark: The last three items of the lemma follow from the following two facts:\\
If $f\in D(A)$ and $t>0$, then $T(t)f-f=\int_0^t T(s)Af\mbox{d}s$.\\
For all $f\in E$ and $t>0$, $\int_0^t T(s)f\mbox{d}s$ is in $D(A)$,
and $A(\int_0^t T(s)f\mbox{d}s)=T(t)f-f$.\\
(The integral is defined as the limit of the Riemann--sums.)\\
Using $\norm{\int f(s)\text{d}s }\leq \int \norm{f(s)} \text{d}s$
the items \ref{itemXY1} and \ref{itemXY2} follow directly
from  $T(t)f-f=\int_0^t T(s)Af\mbox{d}s$.
We get item 
\ref{itemuoehgqtetqw} 
by  defining $f:=\frac{1}{h}\int_0^h T(s) \bar{f} \text{d}s$.

The famous theorem of Hille--Yosida states that uniformly continuous semigroups
are exactly the semigroups of the form $T(t)=e^{t A}$, where $A$ is in
$\Lin(E)$.  If $T(t)=e^{t A}$, then $A$ is the generator of $T$, and $D(A)=E$.
In this case, the constructions in this paper do not result in anything new. So
we will mainly be interested in semigroups that are not uniformly continuous. 

Our basic example for such a semigroup is the following:
\begin{Example}\label{basicexample}
$\mycc$ is the (real) Banachspace of
uniformly continuous, bounded functions from $\bR$ to $\bR$
(with the sup--norm, denoted by $\supnorm{\cdot}$).
$T$ is the translation semigroup, defined by $T(t)f(x):=f(x+t)$.
\end{Example}

Then $T$ is continuous (since each $f\in \mycc$ is uniformly 
continuous), but $T$ is not uniformly continuous: for example,
let $f_k(x)=\sin(k x)$. Then $f_k\in \mycc$ and $\norm{f_k}=1$. 
For all $t>0$ there is a $k$ such that $\norm{T(t)f_k-f_k}>1$,
i.e. $\opnorm{T(t)-\text{Id}_E}>1$. 
$f$ is in $D(A)$ if and only if $f'$ exists and is element of $E$.
In this case, $A(f)=f'$.

\subsection{Robinsonian (nonstandard) analysis:}

We present the concept of nonstandard extensions
following \cite{changkeisler}.

Let $X$ be a set. $P(X)$ denotes the set of all subsets of $X$.
$V_0(X)=X$, $V_{n+1}(X)=V_n(X)\cup P(V_n(X))$,
$V(X)=\bigcup_{n\in\omega} V_n(X)$.
We assume that (a suitable copy of) $\bN$ exists in $X$, 
and that all vector spaces we are interested in are
subsets of $X$.
So all sets of numbers, subspaces of $E$ etc are in $V(X)$.
For technical reasons, we want to treat the elements of $X$
as urelements, i.e. we want to avoid that 
for some $y\in V(X)$ and $x\in X$, $y\in x$.
We can do that by assuming (without loss of generality)
that $X$ is a base--set, i.e. $X$ does not contain the empty set
and $x\cap V(X)$ is empty for all $x\in X$.

A {\em nonstandard universe} is a triple $\la V(X), V(Y), \x:V(X)\fnto V(Y)\ra$
such that:
\begin{itemize}
\item $X$ and $Y$ are infinite base--sets, $X\subset Y$, $\x X=Y$, for all $x\in X$, ${\x x}=x$
\item (non--triviality) for every infinite subset $A$ of $X$: $A\subsetneq {\x A}$
\item (transfer principle) if $\varphi(x_1,\dots,x_n)$
is a $\Sigma_0$ formula and $a_1,\dots,a_n\in V(X)$,
then $V(X)\vDash \varphi(a_1,\dots,a_n)$ iff $V(Y)\vDash \varphi(\x a_1,\dots,
\x a_n)$
\end{itemize}

$\Sigma_0$ formulas are first order formulas $\varphi$ 
(in the Language $\{\in\}$)
such that only quantifiers of the form $\forall x\in y$ and $\exists x\in y$
occur in $\varphi$.

Note that the transfer principle is similar to the well--known
\L o\'{s}' theorem for ultraproduct--constructions.
However, $V(Y)$ cannot be an ultrapower of $V(X)$,
since otherwise \L o\'{s}' theorem would apply to
{\em all} first order formulas, including the formula
``for all natural numbers $n$, there is a decreasing $\in$--chain
of length $n$''. This sentence cannot be true in $(V(Y),\in)$, since
$\in$ is well--founded and  $V(Y)$ has nonstandard natural
numbers. However, $Y$ can be an ultrapower of $X$, and 
this specific kind of nonstandard extension is the most important
one in our context:

Let $\cU$ be a countably incomplete ultrafilter over a set $I$
(i.e. there is a countable family $A_n$ of elements of $\cU$ such
that $\bigcap A_n\notin \cU$). 
Let $Y$ be the ultrapower of $X$ with respect to $\cU$.
(Without loss of generality we can assume that 
$X$ and $Y$ are base--sets).
Then $V(Y)$ can be made a nonstandard extension of $V(X)$ in a way
that the map $\x$ restricted to $X$ is the usual ultrapower
injection, i.e. $\x x$ is the 
equivalence class of the constant function $c_x$.
Non--triviality follows from the fact that $\cU$ is 
countably incomplete, and the transfer principle is 
proved similar to \L o\'{s}' theorem.

Nonstandard extension obtained in this way are called 
{\em bounded ultrapowers}.

In any nonstandard extension there are infinite numbers:
The sets $\xN\supsetneq \bN$ and $\xR\supsetneq \bR$ are called 
the nonstandard natural and real numbers. Sometimes we will 
call $\bN$ and $\bR$ the standard natural (or real) numbers to emphasis the 
difference.

A nonstandard number $r$ is called {\em finite} 
if there is a standard natural number $n$ such that
$\norm{r}<n$ (in $V(\x X)$). So a natural number is finite iff it is
standard.
A nonstandard number that is not finite is called infinite.
A nonstandard real $r$ is called 
{\em infinitesimal} if $r=0$ or $r\neq 0$ and $\frac{1}{r}$ is infinite.

Let $r,s\in{\xR}$. $r\near s$ stands for ``$r-s$ is infinitesimal''.
For each finite $r\in{\xR}$ there is a unique
real number $r'$ such that $r'\near r$. This $r'$ is denoted by $\std(r)$

Similar notation will be introduced for elements of nonstandard normed 
vectorspaces (on page \pageref{defofnearforvectors}).

The following notions are central in nonstandard analysis:\\
$y\in V(Y)$ is called {\em standard} if there is a
$x\in V(X)$ such that $y={\x x}$.\\
$y\in V(Y)$ is called {\em internal} if there is a
$x\in V(X)$ such that $y\in{\x x}$.\\
$y\in V(Y)$ is called {\em external} if $y$ is not internal.

Note that for $y\in Y$, $y$ is standard iff $y\in X$, so the 
notation is compatible with our use of standard natural (or real) numbers.

An example of an internal set that is not standard is the set
of all nonstandard natural numbers less than some $m$, where
$m$ is infinite.

If $A, B$ are internal (or standard),
then so are $A\cup B$, $A\setminus B$, etc.

The set of all infinite natural numbers in not internal.
This is a special case of the so--called 
spillover principle:

\begin{Thm}\label{spillover}{\rm (spillover principle)}
Let $A$ be internal.
\begin{itemize}
\item
$A$ contains arbitrary large finite (i.e. standard) natural numbers iff
$A$ contains arbitrary small infinite natural numbers.
\item
$A$ contains arbitrary small positive standard reals 
iff $A$ contains arbitrary large positive infinitesimals.
\end{itemize}
\end{Thm}

It is important to note that $A$ just has to be internal, it is not
necessary that $A$ is internal.


If $A$ is a set in $V(X)$, then $P^\text{fin}(A)$, the set of all finite subsets of $A$,
is in $V(X)$ as well, and is mapped to $\x P^\text{fin}(A)\in V(Y)$.
If $A$ was infinite, this set will contain more elements than just
the sets $\x B$, where $B\subset A$ finite. 
These new elements are called {\em hyperfinite} (note that 
they do not have to be finite).
So a set $B$ in $V(Y)$ is called {\em hyperfinite}, if there is some 
$A\in V(X)$ such that $B\in {\x P^\text{fin}(A)}$.

The example above, $\{n\in {\x \bN}:\, n<m\}$, is an example of 
an infinite, hyperfinite, internal set.
It is easy to see that all hyperfinite sets are internal.

For the examples \ref{exampleeins} and \ref{examplezwei}, we use
a nonstandard extension such that the standard reals are a 
subset of a hyperfinite set.

More general, we call
a nonstandard extension $V(Y)$ an {\em enlargement} if 
for every $A\in V(X)\setminus X$ there is a hyperfinite $B\in V(Y)$ such that
$\{{\x a}:\, a\in A\}\subset B$.

It is provable (using the Axiom of Choice, of course) that 
for any base--set $X$ there is a index set $I$ and an ultrafilter $\cU$ 
over $I$
such that the bounded ultrapower with respect to $\cU$ is an enlargement.
(For details, see e.g. \cite{changkeisler}. The outline of the 
proof is as follows:\\
1. for every $\kappa$, there are $\kappa$--good filters\\
2. $\kappa$--good filters result in $\kappa$--saturated extensions\\
3. Given $X$, let $\kappa_X$ be the cardinality of $V(X)$.
Then $\kappa_X$--saturated extensions are enlargements.)

The only consequence of the concept of enlargement 
that we are going to use is the following:

\begin{Lemma}\label{lem:factforexampleeins}
  Let $V(\x X)$ be a nonstandard enlargement of the universe.
  Then there is a (infinite) nonstandard natural number
  $k$ such that $\sin(k t)$ is infinitesimal
  for all standard reals $t$.
\end{Lemma}

\begin{proof}
The standard reals $\bR$ are a subset 
of some hyperfinite set $A$ of nonstandard reals.

It is a well--known (classical) fact that for every finite set 
of reals $A$, and every positive $\varepsilon$
there is a natural number $k$ such that for all $t\in A$,
$\sin(k t)<\varepsilon$.
So by transfer of the classical fact, setting 
$\varepsilon$ infinitesimal, we get a 
$k\in\xN$ such that for all $t\in\bR$, $\sin(k t)$ is infinitesimal.
\end{proof}

When we compare nonstandard constructions with classical constructions,
we will assume that our nonstandard
extension is a bounded ultrapower (generated by the same ultrafilter)
--- otherwise it isn't clear what the classical analogon should be.
It should be noted that 
there are nonstandard universes
which are not bounded ultrapowers: An ultrapower is  always 
$\aleph_1$--saturated
(since all filters are $\aleph_1$--good), but the union of a
countable elementary chain of ultrapowers is not 
(see \cite[page 290, 4.4.29]{changkeisler}).
However, every extension can be seen as limit of ultrapowers
(see \cite{changkeisler}), and in
practice, authors  
often restrict their
attention to the case of ultrapowers 
(see e.g.  \cite[page 88]{hurdloeb}).

\section{Ultrapowers and Nonstandard--Hulls}

\subsection{Constructions without Semigroups}


Let $E$ be a normed vectorspace, $I$ an arbitrary set.
We define $\luE$ to be the space of all bounded $E$-valued $I$-sequences with
the sup--norm (denoted by $\supnorm{\cdot}$). Then
$\luE$ is a normed vectorspace. Note that $I$ does not have to be countable,
it can be of any cardinality.

Assume $\cU$ is a filter over $I$ (not necessarily an ultrafilter).
As in \cite{derndinger},
we define $\CU$ to be 
the set of sequences $\fI$ such that for all
$\varepsilon$ there is a set $J\subseteq I$
such that $J\in\cU$ and 
$\norm{f_i}<\varepsilon$ for all $i\in J$.

%

\begin{Lemma}  
$\norm{y+\CU}:=\inf\{ \norm{x} : x\in y+\CU\}$ is a norm on 
$\luED$, and 
$E$ is a Banachspace if and only if $\luE$ is a Banachspace.
If $E$ is a 
Banachspace, then $ \luED$ is a Banachspace as well.
\end{Lemma}

\begin{proof}
If $E$ is a Banachspace, then so is $\luE$ (since 
Cauchy--sequences converge pointwise).
Clearly, \CU\ is a closed subspace of $\luE$.
And any normed (and complete) vectorspace can be factored by a 
closed subspace, resulting in a normed (and complete, resp.) 
vectorspace. 
If $\luE$ is complete, then clearly so is $E$
(otherwise take any non--converging Cauchy--sequence $f_i$ in
$E$, and map it to the sequence $g_i=(f_i,0,0,\dots)$ in
$\luE$).
\end{proof}

Note that $\luED$ can be a Banachspace although $E$ 
(and therefore $\luE$) is not, see the remark after 
theorem \ref{myisom}.

The corresponding nonstandard construction is the following:
Let $V({\x X})$ be a nonstandard universe, and $E \subset X$
a normed vectorspace.
Then $\xE$ is a nonstandard normed vectorspace, and a  
standard vectorspace (without 
canonical norm).

We define the finite part of $\xE$, denoted by $\finxE$, to consist of all
elements $f$ of $\xE$ such that $\norm{f}$ is finite.
\label{defofnearforvectors}$f\near g$ means that 
$\norm{f-g}$ is infinitesimal. 
The infinitesimal part of $\xE$, denoted by $\EO$, consists of
all $f$ such that $\norm{f}$ is infinitesimal, i.e. $f\near 0$.
Clearly, $\EO$ is a (standard) sub--vectorspace of $\finxE$.

$\ED$ is the quotient of $\finxE$ and $\EO$, and the canonical 
quotient map $\finxE\fnto\ED$ is denoted by $\dach{\cdot}$,
i.e. the quotient--class of $f$ is $\fD$, and
$f$ is called representant of $\fD$.

$\ED$ is a normed vectorspace with the norm
$\norm{\fD}:=\std(\norm{f})$, where $f$ is 
any representant of $\fD$ and $\std(r)$ is the
standard part of a finite nonstandard real $r$.
Of course, the vectorspace--operations are defined by 
$\fD+\gD:=\dach{f+g}$ and $\alpha \fD:=\dach{\alpha f}$.  

\begin{Lemma}\label{completehull}
If $V({\x X})$ is a bounded ultrapower,
then $\ED$ is a Banachspace.
\end{Lemma}

A proof can be found e.g. in \cite[page 59]{lindstrom},
noting that a bounded ultrapower is $\aleph_1$--good
and therefore the extension will be $\aleph_1$--saturated.


A straightforward 
calculation proves that for an ultrafilter \cU, the F--product
and the nonstandard hull are the same:


\begin{Thm}\label{myisom} 
Let $\cU$ be an ultrafilter,
$V(\x X)$ the corresponding bounded ultrapower.
Then $\myisomorphism: \luED\fnto \ED$ defined by 
$\myisomorphism(f+\CU)= \dach{f}$ is an isomorphism.
\end{Thm} 

So if $\cU$ is a countably incomplete ultrafilter,
then $\luED$ is a Banachspace (according to lemma \ref{completehull}),
regardless of whether $E$ was a Banachspace or not.

%
%

\subsection{Classical Constructions for Semigroups}

\paragraph{The Maximal Continuous Subspace:}

Assume, $T(t)$ is a semigroup on $E$ (not necessarily continuous),
and assume $M,\omega$ are such that 
$\opnorm{T(t)}\leq M e^{\omega t}$ (for continuous semigroups,
such $M,\omega$ always exist, according to lemma \ref{lem:basicsemigroup}).\\
Define $\MAX{E}:=\{f\in E: \lim_{t\limto 0}\norm{T(t)f-f}=0\}$.

A subspace $F$ of $E$ is called $T$--invariant, if 
$T(t)(f)\in F$
for all $f\in F$ and $t\in \bRp$.

\begin{Lemma}\label{maxsubsp}
$\MAX{E}$ is a closed $T$--invariant subspace of $E$, and
it is maximal in the family of subspaces $F$ of $E$ such that 
$T(t)$ restricted to $F$ is continuous. 
\end{Lemma}

\begin{proof}
\begin{itemize}
\item closed: Assume $f_i, f\in E$ such that 
$f_n\rightarrow f$, and choose an $\varepsilon>0$.
Let $n$ be such that 
$\norm{f_n-f}<\min(\frac{\varepsilon}{6 M},\frac{\varepsilon}{3})$. 
$\norm{T(t)f-f}\leq \norm{T(t)f-T(t)f_n}+\norm{T(t)f_n-f_n}+\norm{f_n-f}$.
If $t<\delta$ (for a suitable $\delta$) then $\opnorm{T(t)}<2 M$
and $\norm{T(t)f_n-f_n}<\frac{\varepsilon}{3}$ (since $f_n\in\MAX{E}$),
so $\norm{T(t)f-f}\leq 2 M \frac{\varepsilon}{6 M}+ 
\frac{\varepsilon}{3}+ \frac{\varepsilon}{3}\leq\varepsilon$.
\item invariant: Assume $f\in\MAX{E},\, s\in\bRpO$. Then
$\norm{T(t)(T(s)(f))-T(s)(f)}=\norm{T(s)(T(t)(f)-f)}\leq
\opnorm{T(s)}\,\norm{T(t)f-f}$.
\item maximal: 
If $T$ is continuous on $F$ and $f\in F$, then
by definition $f$ is an element of $\MAX{E}$. 
\end{itemize}
\end{proof}

Remark: This definition and lemma just isolate a part of the proof that is
given e.g. in \cite{derndinger} for corollary \ref{trivial}. This way we
don't have to repeat the same argument again (e.g. for lemma \ref{hlem44}).

\paragraph{$\mathbf \mT$:}

Let $T(t)$  be continuous, $\opnorm{T(t)}\leq M e^{\omega t}$, and let
$\fI$ be an element of $\luE$. Define
$\TT(t)(\fI):=\la T(t)f_i\ra$.

Since $T(t)$ is bounded, $\TT(t)(\fI)\in\luE$ for all $\fI\in\luE$
and $t\geq 0$.
It is clear that $\TT$ defines a semigroup on $\luE$. 
Also, $\opnorm{\TT(t)}\leq M e^{\omega t}$.

But in general $\TT$ will not be continuous on $\luE$: 
Let $T$ be the translation semigroup on $\mycc$ as in 
\ref{basicexample}, 
$I=\omega$, $\cU$ a free ultrafilter over $\omega$, and
$f_k(x)=sin(k x)\in \mycc$. 
Then $f=\la f_k\ra\in\luE$, but  $\TT(t)f-f$ does not converge.

On can show rather easily that 
$\TT$ is continuous on $\luE$ if and only if $T$ is uniformly continuous on
$E$ (assuming of course that the index set $I$ is infinite).

We define $\mT$ to be the maximal continuous subspace of $\luE$
with respect to $T$. 


\begin{Cor}\label{trivial}
$\mT$ is a closed, $\TT$--invariant subspace of $\luE$.
\end{Cor}

\paragraph{F-Product of a semigroup:}

As we have seen, if $T$ is a semigroup defined on $E$,
under certain assumptions $T$ can be extended to a
(not necessarily continuous) 
semigroup  $\TT$ on $\luE$.
Assume $\fI\in\CU$ (i.e.
for all $\varepsilon\in\bRp$, 
$\norm{f_i}<\varepsilon$ holds on a filter--set).
$\TT(t)(\fI)\in\CU$ for all $t\geq 0$, since
$\TT(t)$ is a bounded operator.
This means that $\CU$ is $\TT$--invariant, and that
we can define $\TT$ on $\luED$.
Of course, $\TT$ is a semigroup on $\luED$ 
(and generally not continuous).

$\CU$ is not a subspace of $\mT$.
(For example, any $\fI$ such that 
$f_i=0$ on a filter--set is in $\CU$,
but not necessarily in $\mT$.)
However, it is easy to see that 
$\CU\cap\mT$ is a closed subset of $\mT$, and
$\mTD$ is a subspace of $\luED$.

To be more exact:

\begin{Lemma}\label{mTDsubsetluED}
$\myisometry: \mTD\fnto\luED$ 
defined by 
$\myisometry(\fI+\CU\cap\mT)=\fI+\CU$ is a (well--defined) 
injective isometry.
\end{Lemma}

\begin{proof}
See \cite[page 163]{krupa}.
\end{proof}

\begin{Lemma}
$\mTD$ is a closed, $\TT$--invariant subspace of $\luED$.\\
$\TT$ restricted to $\mTD$ is continuous.
\end{Lemma}

\begin{proof}
The quotient $\mTD$ of the Banachspace $\mT$ is a 
Banachspace, and so is its isometric image. The rest is clear.
\end{proof}

In general, $\mTD$ is {\em not} the maximal continuous subspace 
of $\luED$. 
This is shown in example \ref{exampleeins} (using 
that $\mTD$ corresponds to $\ETD$, see lemma \ref{lem:hcor66}).

\subsection{Nonstandard Constructions for Semigroups}

Let $V(\x X)$  be a nonstandard universe, $E\subset X$
a Banachspace.

Assume $A:E\fnto E$ is a continuous linear operator, and
$f$ an element of $\finxE$, $g$ an element of 
$\EO$ (the finite and infinitesimal parts of
the nonstandard vector space $\xE$, resp.). 
Then ${\x A}f$ (and ${\x A}g$) are again finite (and
infinitesimal, resp.).
Therefore $A$ defines a bounded linear operator
$\dach{A}: \ED\fnto\ED$
by $\dach{A}(\fD)=\dach{A f}$.
(If $A$ is not bounded, then generally $\dach{A}$ cannot be 
defined in a canonical way.)

If $T$ is a semigroup, then for all positive reals $t$, 
$T(t)$ is continuous, 
and $\dach{T}(t):=\dach{T(t)}$ is a 
semigroup on $\ED$.

If $\opnorm{T(t)}<M e^{\omega t}$ (this is always the case
if $T$ is continuous), then 
an alternative way to define $\dach{T}$ is the following:
$\xT$ is a nonstandard semigroup on $\xE$.
If $t$ is a positive, finite nonstandard real, then
$\opnorm{\xT(t)}<M e^{\omega t}$, so we get:
If $f\in \xE$ is finite (or infinitesimal), 
then so is $\xT(t)f$. Therefore $\dach{T}(t):\ED\fnto\ED$
is well--defined
by $\TD(t)(\fD)=\dach{\xT(t)(f)}$, even for finite nonstandard $t$.
If $t$ is standard, then the two definitions are equivalent.


Assume $T$ is continuous, i.e. for all $f$ in $E$ and all
$\varepsilon>0$ there is a $\delta > 0$ such that 
$\norm{T(t)f-f}<\varepsilon$
for all $t<\delta$.
We apply the transfer principle to this sentence, and get:
For all $f\in \xE$ and all $\varepsilon\in \xRp$, 
there is a $\delta\in \xRp$ such that 
$\norm{\xT(t)f-f}<\varepsilon$
for all positive reals $t<\delta$: 
However, if $\varepsilon>0$ is a standard real, then
it is not guaranteed that $\delta>0$ can be chosen to be
standard as well ($\delta$ could 
be infinitesimal). 
We define
two subspaces of $\finxE$ (the finite part of $\xE$):
$\ET$ consists of all the vectors $f$ with the following property:
If $\varepsilon>0$ is standard, then there is a standard $\delta>0$ 
such that $\norm{\xT(t)f-f}<\varepsilon$ 
for all nonstandard $t<\delta$.
$\Emax$ consists of all $f$ satisfying the same condition for
standard $t<\delta$ only. I.e. a finite $f\in \xE$ is in 
$\Emax$ iff for all $\epsilon\in\bRp$ there is a $\delta\in\bRp$
such that 
$\norm{\xT(t)f-f}<\varepsilon$
for all positive standard reals $t<\delta$.

It is easy to see that a finite $f\in \xE$ is in $\ET$ iff
for all positive infinitesimal t, $\norm{\xT(t)f-f}$ is infinitesimal.


Clearly, $\EO$ is a subspace of $\ET$, which is in turn a subspace 
of $\Emax$.
So, if $f\in \ET$ (or $\Emax$) and $f\near g$,
then $g\in \ET$ (or $\Emax$, resp.), and we can define the
quotient spaces 
$\ETD=\myfrac{\ET}{\EO}$ and $\EmaxD=\myfrac{\Emax}{\EO}$.
Remember that (other then $\finxE$), $\ET$ is a (canonically)
normed vectorspace.

\begin{Lemma}\label{hlem44}
\begin{enumerate}
\item $\EmaxD$ is the maximal continuous subspace of $\ED$ with respect to 
$\TD$.
\item $\EmaxD$ is a closed, $\TD$--invariant subspace of $\ED$
\item $\ETD$ is a closed, $\TD$--invariant subspace of $\EmaxD$.
\end{enumerate}
\end{Lemma}

\begin{proof}
\begin{enumerate}
\item Assume, $\lim_{t\limto 0}\norm{\TD(t)\fD-\fD}=0$,
and fix $\varepsilon\in\bRp$. Then there 
is a $\delta\in\bRp$
such that for all positive reals $t<\delta$, 
$\norm{\TD(t)\fD-\fD}<\frac{\varepsilon}{2}$.
Assume that $f$ is a representant of the quotient $\fD$, and
that $t<\delta$. Then 
$\std(\norm{\xT(t)f-f})<\frac{\varepsilon}{2}$, and therefore
$\norm{\xT(t)f-f}<\varepsilon$ in $\xR$, 
so $f\in\Emax$ and $\fD\in\EmaxD$.
\item this follows from 1. and lemma \ref{maxsubsp}.
\item Assume that $\fD_n\in\ETD$, and $\fD_n\limto \fD$ in $\ET$.
Let $f_n$ and $f$ be representants of $\fD_n$ and $\fD$, resp.
Assume $h$ is a positive infinitesimal. 
Then $\norm{T(h)f-f}\leq\norm{T(h)f-T(h)f_n}+\norm{T(h)f_n-f_n}+\norm{f_n-f}$,
which is smaller than every positive standard $\epsilon$, since
$\norm{f-f_n}$ gets arbitrary small,
$\norm{T(h)f-f}\leq \opnorm{T(h)}\, \norm{f-f_n}$, and
$\norm{T(h)f_n-f_n}$ is infinitesimal.

To show the invariance, assume that 
$\fD\in \ETD$, and let
$f\in\ET$ be a representant. Let $t$ be a positive standard real,
and $h$ is a positive infinitesimal.
$\norm{T(h)T(t)f-T(t)f}\leq
\opnorm{T(t)}\, \norm{T(h)f-f}$, which is infinitesimal.
So $T(t)f\in\ET$, and therefore $\TD(t)\fD\in\ETD$.

\end{enumerate}
\end{proof}

Remark: In \cite{wolff_nsfa}, an alternative way to construct \ET\ is
presented. Theorem 4.6.1 there corresponds to our \ref{wid}.

In general,  $\ETD$ is a proper subset of $\EmaxD$
(or equivalently: $\ET$ is a proper subset of $\Emax$).
To see this, we bring
an example similar to \cite[page 165]{krupa} or \cite[page
209]{wolff_spectral}:

\begin{Example}\label{exampleeins}
  Let $T$ be the translation semigroup on $\mycc$ (as in example 
  \ref{basicexample}), 
  $V(\x X)$ a enlargement, $k$ such that
  $\sin(k t)$ is infinitesimal
  for all $t\in \bR$ (see lemma \ref{lem:factforexampleeins}).
  Define the element $f$ of $\x\mycc$ by  $f(x)=\sin(k x)$.
  Then $\norm{f}=1$, 
  and $f$ is an element of $\Emax$, but not of $\ET$.
\end{Example}

\begin{proof}
To see that $f\in \Emax$, it is enough to show that 
for all positive (standard) reals $t$,
  $\norm{T(t)f-f}$ is infinitesimal.
For any standard reals $t$,
$\sin(k(x+t))-\sin(k x)=
\sin(k x)(\cos(k t)-1) +\sin(k t)\cos(k x)$.
$\sin(k x)$ and $\cos(k x)$ are finite, and
$\sin(k t)$ and $(\cos(k t)-1)=(\sqrt{1-\sin^2(k t)} -1)$  are infinitesimal.
Therefore $\norm{T(t)f-f}=\supnorm{\sin(k(x+t))-\sin(k x)}$ is infinitesimal.

To see that $f\notin \ET$, it is enough to note that
$h:=\frac{1}{k}$ is infinitesimal, but 
$\norm{T(h)f-f}=\supnorm{sin(k(x+h))-sin(k x)}=
\supnorm{sin(x+1)-sin(x)}$ is not.
\end{proof}

\subsection{Relation of Classical and Nonstandard Constructions}

Assume $V(\x X)$ is the bounded ultrapower of an ultrafilter
$\cU$ over $I$. Then the 
isomorphism $\myisomorphism$ of theorem \ref{myisom}
allows us to identify  $\luED$ and $\ED$. 
We can also consider \mTD\ to be a subspace of $\luED$ (via
the injective isometry \myisometry\ of lemma \ref{mTDsubsetluED}).

\begin{Lemma}
$\EmaxD$ is the maximal continuous subspace of $\luED$ 
with respect to $\TT$ 
\end{Lemma}

(More formally, one would have to write 
$\EmaxD=\myisomorphism(\MAX[\TT]{\luED})$.)

\begin{proof}
This follows from lemma \ref{hlem44} ,
and the fact that 
$\myisomorphism(\TT(t)f) =\TD(t)(\myisomorphism(f))$
for all $f\in\luED$.
\end{proof}

\begin{Lemma}\label{lem:hcor66}
$\mTD\subset\ETD$
\end{Lemma}

(Again, more formally this should be written as
$\myisomorphism(\myisometry(\mTD))\subset\ETD$.)

\begin{proof}
%
%
%
%
%

By definition, $\fI\in \mT$  
iff for all reals $\epsilon>0$ 
there is a real $\delta>0$ such that $\supnorm{T(t)f_i-f_i}<\epsilon$
for all reals $t<\delta$.
Let $\fD\in \ED$ correspond to $\fI$.
Assume that  $h$ is an infinitesimal  nonstandard real. We want to 
show that for an arbitrary fixed positive real $\epsilon$,
$\norm{\xT(h)\fD-\fD}<\epsilon$. Let $\delta$ be the 
real corresponding to $\epsilon$ as above.
$h$ corresponds to a sequence of reals $\hI$ ($i\in I$)
such that $U=\{i\in I:\, h_i<\delta\}$
is in the ultrafilter.
$\xT(h)\fD-\fD$ corresponds to the sequence $T(h_i)f_i-f_i$.
If $i$ is in $U$, then $\norm{T(h_i)f_i-f_i}<\epsilon$,
and $U$ is a filter--set, therefore $\norm{\xT(h)\fD-\fD}<\epsilon$.
\end{proof}

It is not immediately clear whether $\mTD$ is always 
{\em identical} to $\ETD$.

One could construct a counterexample under the following assumption:

%

\begin{Assumption}\label{assumption}
Fix $\omega, M, \eta, \varepsilon_1,\varepsilon_2,\dots, \delta_1,\delta_2,\dots$.
Assume that for all $n,m\in\bNp$ there is a Banachspace $E^n_m$, 
a continuous 
semigroup $T^n_m$ with $T^n_m(t)<M e^{\omega t}$ 
and a $f^n_m\in E^n_m$ such that 
\begin{itemize}
\item $\norm{f^n_m}=1$
\item 
        $\norm{T^n_m(t)f^n_m-f^n_m}<\frac{1}{n}$
for all $1\leq i\leq n$, $t<\delta_i$ 
\item For all $\bar{f}\in E^n_m$ such that $\norm{\bar{f}-f^n_m}<\eta$
exists a $t<\frac{1}{m}$ such that $\norm{T^n_m(t)\bar{f}-\bar{f}}>\varepsilon_n$
\end{itemize}
\end{Assumption}

If this assumption holds, then let the index set $I$ be 
$\bNp\times \bNp$, and define
$\bar{E}=\lunendl(E^i_j)$ with the sup--norm.
(I.e. $\la f_{i,j}\ra \in \bar{E}$ iff
$f_{i,j}\in E^i_j$ and
$\norm{f_{i,j}}$ is bounded.)
$T(t)\la g_{i,j}\ra:=\la T^i_j(t)g_{i,j}\ra$ is a semigroup on $\bar{E}$.
Let $E$ be the maximal continuous subspace of $\bar{E}$,
and define $h^n_m=\la h^n_m(i,j) \ra \in E$ by\\
$h^n_m(i,j)= 
\begin{cases}f^n_m & 
\text{if $i=n, j=m$,}\\ 
0  & 
\text{otherwise}.
\end{cases}$.
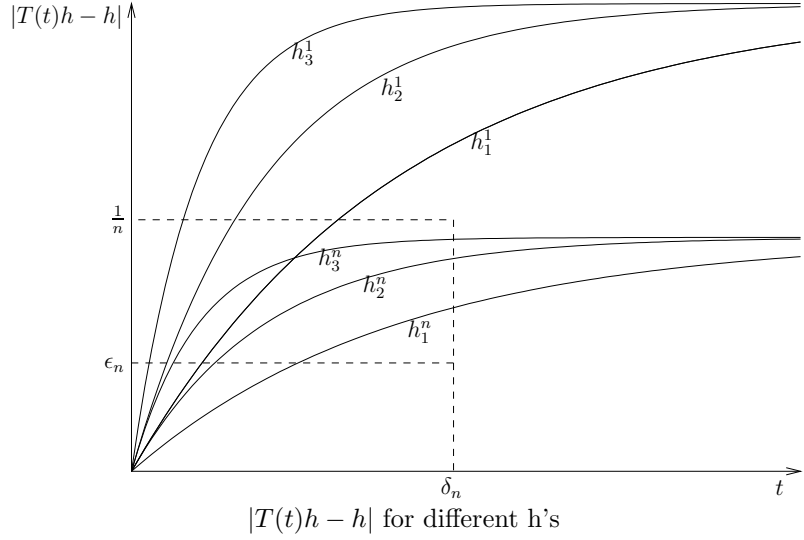
\begin{figure}[h]
    \centering
    \scalebox{0.75}{\input{counterexample.tex}}
    {$\norm{T(t)h-h}$ for different h's}
\end{figure}

Assume that $\cU$ is an ultrafilter over $I=\bNp \times \bNp $, 
and that there is a partition
$\{\sigma_{i,j}\}$ $(i,j\in\bN)$ of $I$ such that for all $n\in\bN$ 
and all functions 
$f:\bN\fnto\bN$, $\bigcup_{i>n,\,j>f(i)} \sigma_{i,j}\in\cU$
(it is easy to see that such a countably incomplete filter exists,
e.g. apply a suitable filter--basis to \cite[theorem A.4]{lindstrom}).

Then define $k=\la k_{i,j}\ra\in\bar{E}$
by $k_{i,j}=h^n_m$, where $n,m$ is the
(unique) pair of natural numbers such that 
$(i,j)$ is an element of $\sigma_{n,m}$.

Now consider the nonstandard extension defined by $\cU$,
and interpret $k$ as nonstandard element of $\xE$.

Clearly $\norm{k}=1$, since
$\norm{k_{i,j}}=1$ for all $(i,j)$.
For all natural numbers $n$, $0<t<\varepsilon_n$ and $(i,j)\in 
\bigcup_{r>n,\,s\in\omega}\sigma_{r,s}$:
$\norm{T(t)k_{i,j}-k_{i,j}}<\frac{1}{n}$. 
So if $\la r_{i,j}\ra$ is infinitesimal, then 
so is $\norm{\xT(r)k-k}$, i.e. $k\in\ET$.

Assume that the equivalence class of $k$ is in $\mTD$.
Pick a representant $\la\bar{k}_{i,j}\ra\in \mT$,
and $U\in \cU$ such that 
$\norm{\bar{k}_{i,j}-k_{i,j}}<\eta$
for all $(i,j)\in U$.
Since $U\in\cU$, there is an $n_0$ such that 
for all $m$ there is a $m'$ such that
$\sigma_{n_0,m'}\cap U\neq\emptyset$. 
Since $\la \bar{k}_i\ra\in\mT$, there is a $m_0\in\bN$ such that 
$\norm{T(t)\bar{k}_i-\bar{k}_i}
<\delta_{n_0}$
for all $(i,j)\in I$ and $0<t<\frac{1}{m_0}$.
If ${m_0}'>m_0$ and $(i_0,j_0)\in \sigma_{n_0,{m_0}'}\cap U$, 
then
there is a $t\in\bRp$ such that $t<\frac{1}{{m_0}'}<\frac{1}{m_0}$ and 
$\norm{T(t)\bar{k}_i-\bar{k}_i}\geq\delta_{n_0}$, a contradiction.

However, our  investigation of  the nonstandard generator of a semigroup
will show that $\mTD=\ETD$.

\section{The Generator}

\subsection{The Generator in $\xE$}

If $T$ is a continuous semigroup, then 
$\xT(t)$ maps the finite (or infinitesimal) elements 
of $\xE$ to finite (or infinitesimal, resp.) elements
for all finite $t\geq 0$.
That is not true for the (unbounded) linear map
$A$, as the following trivial example shows:

\newcommand{\INF}{k}
\begin{Example}\label{exampletrivial}
  Let $T$ be the translation semigroup on $\mycc$ (as in example
  \ref{basicexample}),
  $V(\x X)$ a nonstandard extension, $\INF$ an infinite natural
  number. 
  Then there is an infinitesimal $f_1\in\xE$ such that $f_1\in \xD(A)$ and 
  $\xA(f)$ is infinite.
\end{Example}
\begin{figure}[h]
    \centering
    \scalebox{0.75}{\input{example1.tex}}\\
    {Example \ref{exampletrivial}: $f_1$ and $f_1'=\xA(f_1)$}
    \label{fig:example1}
\end{figure}

\begin{proof} Define $f_1$ by case distinction:\\
$f_1=\begin{cases}
0                                 &           x<0  \\
\INF^3 x^2                        &  0   \leq x<\frac{1}{\INF^2 } \\
-\INF^3 x^2+4\INF x-\frac{2}{\INF}& \frac{1}{\INF^2}\leq x<\frac{2}{\INF^2} \\
\frac{2}{\INF}                    &  \frac{2}{\INF^2}\leq x
\end{cases}$.

Then
$\xA f_1=\begin{cases}
0                   &                       x<0  \\
2\INF^3 x           &  0               \leq x<\frac{1}{\INF^2 } \\
-2\INF^3 x+4\INF    &  \frac{1}{\INF^2}\leq x<\frac{2}{\INF^2 } \\
0                   &  \frac{2}{\INF^2}\leq x
\end{cases}$.
\end{proof}

Let $V(\x X)$ be a nonstandard universe.
Then the subspace $D(A)$ of $E$ is mapped to a subspace $\x D(A)$ of
$\xE$.

\begin{Lemma}\label{LemET}
Assume, $f$ and $g$ are elements of $\x D(A)$
\begin{enumerate}
\item If $f$ and $\xA (f)$ are both finite, then $f\in\ET$.
\item\label{itemweqtqw} If $f$ is finite and $\xA (f)\in\ET$, then 
      $\frac{\xT(h)f-f}{h}\near\xA(f)$
      for all positive infinitesimal reals $h$.
\item If $k\in \xE$ and  
      $\frac{\xT(h)f-f}{h}\near k$ 
      for all positive infinitesimal reals $h$,
      then $k\near\xA(f)$
\item If $f$, $g$, $\xA(f)$ and $\xA(g)$ all are in $\ET$, and $f\near g$, 
      then $\xA(f)\near\xA(g)$
\end{enumerate}
\end{Lemma}

\begin{proof}
\begin{enumerate}

\item By transfer of \ref{lem:basicsemigroup}.\ref{itemXY1} we get:
$\norm{\xT(t)f-f}\leq 
t\cdot \sup_{s\leq t}\opnorm{\xT(s)}\cdot\norm{\xA(f)}\leq t M (e^{\omega t}+1)\norm{\xA(f)}$,
which is infinitesimal for infinitesimal $t$.

\item 
For every standard $\varepsilon>0$ and infinitesimal $s$, we have 
$\norm{\xT(s)\xA(f)-\xA(f)}<\varepsilon$.
Therefore 
transfer of \ref{lem:basicsemigroup}.\ref{itemXY2}
implies that
$\norm{\frac{\xT(h)f-f}{h}-\xA(f)}\leq \varepsilon$ for infinitesimal $h$.
Since $\varepsilon$ was arbitrary, $\norm{\frac{\xT(h)f-f}{h}-\xA(f)}$
is infinitesimal.


\item 
By transfer of the definitions of $A$ and $D(A)$, we get the following:
For all positive nonstandard reals $\varepsilon$
there is a positive nonstandard real $\delta$ such that
$\norm{\frac{\xT(h)f-f}{h}-\xA f}< \epsilon$ for all $h\leq \delta$.
Let $\varepsilon$ be a fixed positive infinitesimal, and 
choose an infinitesimal $h\leq \delta$
for the appropriate $\delta$. Then
$\frac{\xT(h)f-f}{h}\near \xA f$, and by assumption 
$\frac{\xT(h)f-f}{h}\near k$, therefore $k\near \xA$.


\item 
Since $\xA f\in \ET$, $\frac{\xT(h)f-f}{h}\near\xA f$ for all
infinitesimal $h$, and the same applies to $g$ and $\xA g$.
So we have 
$\norm{\xA f-\xA g}\leq
\norm{\xA f-\frac{\xT(t)f-f}{t}}+
\norm{\frac{\xT(t)f-f}{t}-\frac{\xT(t)g-g}{t}}+
\norm{\frac{\xT(t)g-g}{t}-\xA g}\leq\varepsilon$.
So we just have to show that
$\norm{\frac{\xT(t)f-f}{t}-\frac{\xT(t)g-g}{t}}
=\norm{\frac{\xT(t)(f-g)-(f-g)}{h}}$ is infinitesimal for infinitesimal $h$.
Transfer of \ref{lem:basicsemigroup}.\ref{itemXY1} shows that
$\norm{\frac{\xT(t)(k)-(k)}{h}}\leq 
\sup_{s\leq h}\opnorm{\xT(s)}\, \norm{\xA(k)}$
for all positive $h$ and $k\in \xD(A)$.
Now  apply this to $k=f-g$.

\end{enumerate}
\end{proof}


As we have already seen in example \ref{exampletrivial}, 
$f\near g$ does not imply $\xA(f)\near\xA(g)$
(set $f:=f_1$, $g:=0$). 
Also, it is not enough to assume that
$f$, $g$, $\xA(f)$ and $\xA(g)$ are all 
in $\Emax\cap \xDA$ for the implication to hold:

\begin{Example}\label{examplezwei}
Assume $T$, $k$ and $f$ are as in example \ref{exampleeins}.
Define $g(x):=-\frac{\cos{k x}}{k}$. Then $g\near 0$. Clearly
$g\in \Emax$ (it is even in $\ET$),
and $\xA(g)=f$, i.e. $\xA(g)\in \Emax$, but $\xA(g)$ is not 
infinitesimal.
\end{Example}

The same example shows that
$g\in \ET$ does not imply $\xA(g)\in \ET$.

\subsection{The Generator in $\ETD$}

As we have already seen in lemma \ref{hlem44},
$\TD$ is a continuous semigroup on $\ETD$ (since $\ETD$
is a subspace of $\EmaxD$).
The generator is called $\AD$.

\begin{Thm}\label{wid}
$\fD$ is element of $D(\AD)$, iff there is a representant
$f$ of $\fD$ such that $f$ is in $D\xA$ and $\xA f$ is in $\ET$.
In this case, $\AD(\fD)=\dach{\xA(f)}$.
\end{Thm}

\begin{proof}
Assume $\fD$ has a representant $f$ as in the lemma.
We have to show that $\fD\in D(\AD)$ and $\AD(\fD)=\dach{\xA(f)}$, 
i.e. we fix a (standard) real $\varepsilon>0$ and have to find 
a $\delta>0$ such that $\norm{\frac{\TD(h)\fD-\fD}{h}-\dach{\xA(f)}}<
\varepsilon$ for all $h<\delta$. 
By lemma \ref{LemET}.\ref{itemweqtqw},
$n(h):=\norm{\frac{\xT(h)f-f}{h}-\xA(f)}$ is infinitesimal for
infinitesimal $h$. 
Note that $n$ is an internal function, since it is element of 
$\x\{h \mapsto \norm{\frac{1}{h}(T(h)f-f)-A(f)}:\, f\in D(A)\}$.
Therefore the set $X=\{h:\, n(h)\geq \varepsilon\}$ is internal as well,
and 
we can apply the spillover principle \ref{spillover}:
Assume toward a contradiction that for all standard 
$\delta>0$ there was a standard $h<\delta$ such that 
$n(h)\geq \varepsilon$, i.e. there are arbitrary small standard
reals in $X$. Then there is a infinitesimal real in $X$ as well,
a contradiction. Therefore $\fD\in D(\AD)$.


Assume on the other hand that $\fD$, $\gD$ are elements of $\ETD$
such that $\AD(\fD)=\gD$.
We want to show that $f\in D\xA$ and $\xA f\in\ET $
for some representant $f$ of $\fD$.
Let $\bar{f}, \bar{g}$ be arbitrary representants of $f$ and $g$, resp.
(According to example \ref{exampletrivial},
we cannot hope that $\bar{f}$ already has the required properties.)
For all nonstandard natural numbers $m$ 
define $c_m$ to be the set if nonstandard natural numbers $n>m$
such that $\norm{ \frac{\xT(\frac{1}{n})f'-f'}{\frac{1}{n}}-g'}<\frac{1}{m}$.
For all standard natural numbers $m$, the set $c_m$ is not empty, 
since $\AD(\fD)=\gD$.
The set $X=\{m:\, c_M\neq\emptyset\}$ is internal, and contains
all standard natural number, therefore it contains an infinite natural
number as well (according to the spillover principle).
Let $m$ be such an infinite number, and $n$ an element of $c_m$.
Now we apply the transfer of \ref{lem:basicsemigroup}.\ref{itemuoehgqtetqw},
setting $h:=\frac{1}{n}$.  
So we get a $f$ in $D(\xA)$ such that $\norm{\bar{f}- f}\leq
\sup_{s\leq h}(\norm{\xT(s)\bar{f}-\bar{f}})$, which is 
infinitesimal, since $\bar{f}\in \ET$ and $s$ is infinitesimal. 
Therefore $f$ is a representant of $\fD$ as well.
Also, $\xA(f)=\frac{\xT(h)\bar{f}-\bar{f}}{h}$, so
$\norm{\xA(f)-\xA(\bar{g})}\leq \frac{1}{m}$
is infinitesimal, therefore $\xA(f)\in \ET$ as well.

\end{proof}

Let $V(\x X)$ be a bounded ultrapower, $T$ continuous 
with generator $A$. 

\subsection{The equivalence of F--product and nonstandard hull}

\begin{Lemma}\label{lem:ScETDeqmTD}
	\ETD=\mTD
\end{Lemma}

\begin{proof}

Assume $\mTD\subsetneq\ETD$. Since $\mTD$ is a closed subspace of
$\ETD$ and $D(\AD)$ dense in $\ETD$, there must be a 
$\fD$ which is element of $D(\AD)$ but not $\mTD$.
We chose a representant $f$ of $\fD$ 
such that $f\in D(\xA)$, $\xA f\in \ET$.
Since $f$ is an element of the ultraproduct, $f$ is 
an equivalence class of a sequence $f_i$ ($i\in I$) such 
that for all $i$
$f_i\in D(A)$, $\la A(f_i)\ra \in \luE $. 
But then 
$\norm{\TT(h)f-f}=\sup_{i\in I}(\norm{T(h)f_i-f_i})\leq
\leq h\cdot  \sup_{i\in I}(\norm{A(f_i)})\cdot \sup_{s\leq h}(\opnorm{T(s)})$
(according to \ref{lem:basicsemigroup}.\ref{itemXY1}), so $f\in \mT$,
a contradiction.
\end{proof}

Combining this with \ref{lem:hcor66} and \ref{hlem44}, using 
the isometry $\myisomorphism$ of lemma \ref{myisom}, we get:
\begin{Thm}\label{main}
$\ED$ is isomorphic to $\luED$, and $\ETD$ is isomorphic to
$\mTD$.
\end{Thm}

So we get the following picture
($p$ are the canonical projections from a space to its quotient,
and $\dach{\cdot}$ maps a finite element $f$ of $\xE$ to 
the equivalence class $\fD$. 
The maps labeled with $p$ and $\dach{\cdot}$ are 
surjective, the ones labeled with $\cong$ isometries.
$>$ denotes the subspace--relation, 
for $>^1$ the isometry $\myisometry$ of lemma \ref{mTDsubsetluED} is
used):

\newcommand{\Md}{\ \dach{ }}
\newcommand{\Mi}{\myisomorphism}

\xymatrix{
\luE   \ar[d]^p               &  & >   
        &    & \mT \ar[d]^p   \\
\luED  \ar@{<->}[d]^\Mi_\cong & > &\MAX{\luED}  \ar@{<->}[d]^\Mi_\cong 
        & >^1 & \mTD  \ar@{<->}[d]^\Mi_\cong   \\
\ED                           & > &\EmaxD                         
        & >   & \ETD           \\
\finxE \ar[u]_\Md             & > &\Emax               \ar[u]_\Md 
        & >   & \ET \ar[u]_\Md
}


(Remark: The entries of the first three rows are Banachspaces, the 
last row consists of vectorspaces.)

So we know that assumption \ref{assumption} cannot hold. 
This 
proves 
\begin{Cor}\label{little}
Assume that $T^n_m$ are contraction semigroups on $E$
(i.e. for all $t\in \bRp$, $\opnorm{T(t)}\leq 1$), and
and $f^n_m\in E$ are such that
$\norm{f^n_m}=1$ and 
        $\norm{T^n_m(t)f^n_m-f^n_m}<\frac{1}{n}$
for all $t<\frac{1}{n}$.
Then there is a 
$f'\in E$, and $n, m$, s.t. $\norm{f'-f^n_m}<\frac{1}{2}$ and 
      $\norm{T^n_m(t)f'-f'}<\frac{1}{n}$
for all $t<\frac{1}{m}$.
\end{Cor}


%
%
%






\bibliographystyle{amsalpha}
\bibliography{fproducts}
\end{document}